\documentclass[12pt,a4paper,twoside]{article}

\usepackage[centertags]{amsmath}
\usepackage[T1]{fontenc}
\usepackage{amsfonts}
\usepackage{amssymb}
\usepackage{ntheorem}
\usepackage{ulem}
\usepackage{epsfig}
\usepackage{subfigure}
\usepackage[all]{xy}
\usepackage{changebar}
\usepackage{color}
\usepackage{mathrsfs}
\usepackage[boldsans]{ccfonts}
\usepackage{array}
\usepackage{tabulary}
\usepackage{supertabular}
\usepackage{calc}
\usepackage{eepic}
\usepackage{titlesec}
\usepackage{graphicx}
\usepackage{enumitem}

\theoremstyle{break} \newtheorem{theorem}{Theorem}[section]
\theoremstyle{break} 
\theoremstyle{break} 
\theoremstyle{nonumberbreak}  
\theoremstyle{break} \newtheorem{lemma}[theorem]{Lemma}
\theoremstyle{break} \newtheorem{corollary}[theorem]{Corollary}
\theoremstyle{break} 
\theoremstyle{break} 
{\theorembodyfont{\rmfamily}
\newtheorem{remark}[theorem]{Remark}
{\theorembodyfont{\rmfamily}}
\theoremstyle{break} \newtheorem{problem}[theorem]{Problem}
\theoremstyle{break} 
\theoremstyle{break} \newtheorem{proposition}[theorem]{Proposition}
\theoremstyle{break} 
\numberwithin{equation}{section}

\newcommand{\hide}[1]{}
\newcommand{\N}{{\mathbb{N}}}
\newcommand{\R}{{\mathbb{R}}}
\newcommand{\D}{{\mathbb{D}}}
\newcommand{\C}{{\mathbb{C}}}

\def\Re{\mathop{{\rm Re}}}

\begin{document}

\begin{center}
{\Large \bf Maximal Blaschke products}
\footnotetext{2000 Mathematics Subject
 Classification: 30H05, 30J10, 35J60, 30H20, 30F45, 53A30\\ This research has been
supported by the Deutsche Forschungsgemeinschaft (Grants: Ro 3462/3--1 and Ro
3462/3--2) \hfill{Date: January 31, 2013}}

\end{center}
\renewcommand{\thefootnote}{\arabic{footnote}}
\setcounter{footnote}{0}

\smallskip

\begin{center}
{\large Daniela Kraus and Oliver Roth}\\[3mm]
\end{center}

\smallskip

\begin{center}
\begin{minipage}{13.2cm}
{\bf Abstract.}{ \small
We consider the classical problem of maximizing the derivative at a fixed point over the
set of all bounded analytic
functions in the unit disk with prescribed critical points. We show that the
extremal function is essentially unique and always an indestructible Blaschke product. This
result extends the Nehari--Schwarz Lemma and leads to a new class of
Blaschke products called maximal Blaschke products. We establish a number of
properties of maximal Blaschke products, which indicate that maximal Blaschke
products constitute an appropriate infinite generalization of the class of finite Blaschke products.
}
 \end{minipage}
\end{center}

\medskip

\section{Introduction and Results}

Let $H^{\infty}$ denote the space of all functions analytic and bounded
in the 
unit disk  $\D:=\{z \in \C \, : \, |z|<1\}$ equipped with the norm 
$$ {||f||}_{\infty}:=\sup \limits_{z \in \D} |f(z)|<\infty \, , \qquad f \in
H^{\infty} \, .$$
A sequence  $\mathcal{C}=(z_j)$ in $\D$ is called an $H^{\infty}$ critical set, 
if there exists a  nonconstant function $f$ in $H^{\infty}$ whose critical points are
precisely the points on the sequence $\mathcal{C}$ counting
multiplicities. This means that
 if the point $z_j$ occurs $m$ times in the sequence, then $f'$ has a zero at $z_j$ of
precise order $m$, and $f'(z)\not=0$ for every $z \in \D \backslash \mathcal{C}$.
For an $H^{\infty}$ critical set $\mathcal{C}$ we define 
 the subspace
$$ \mathcal{F}_{\mathcal{C}}:=\left\{f \in H^{\infty} \, : \, f'(z)=0 \text{ for
  any } z \in \mathcal{C} \right\} $$
of all functions $f \in H^{\infty}$ such that any point of the sequence
$\mathcal{C}$ is a critical point of $f$  (with at least the prescribed multiplicity).

\medskip

Our first theorem shows that the set $\mathcal{F}_C$ always contains a Blaschke product
whose critical set is precisely the sequence $\mathcal{C}$. In fact, more
is true:

\begin{theorem} \label{thm:2}
Let $\mathcal{C}=(z_j)$ be an $H^{\infty}$ critical set and let
$N$ denote the number  of times that $0$ appears in the sequence ${\cal  C}$.
Then the  extremal problem
\begin{equation} \label{eq:ex1}
 \max \big\{ \Re f^{(N+1)}(0) \, : f \in {\cal F}_{\cal C} , \,
 {||f||}_{\infty} \le 1 \big \}
\end{equation}
has a unique solution $B_{\mathcal{C}} \in \mathcal{F}_{\mathcal{C}}$.
The extremal function $B_{\mathcal{C}}$ is an indestructible Blaschke product  with
critical set $\mathcal{C}$ and is normalized by
$B_{\mathcal{C}}(0)=0$ and $B^{(N+1)}_{\mathcal{C}}(0)>0$. If $\mathcal{C}$
is a finite sequence consisting of $m$ points, then $B_{\mathcal{C}}$ is a finite
Blaschke product of degree $m+1$; otherwise, $B_{\mathcal{C}}$ is an infinite
Blaschke product.
\end{theorem}

Note that Theorem \ref{thm:2} says that the critical points of the extremal function
$B_{\mathcal{C}}$ are \textit{exactly} the points of the sequence
$\mathcal{C}$ with prescribed multiplicity, so there are no ``extra critical
points'' and $\mathcal{C}$ is the critical set of $B_{\mathcal{C}}$.

\medskip

The crucial part of Theorem \ref{thm:2} is the assertion that the extremal
function $ B_{\mathcal{C}}$ is always 
 an \textit{indestructible Blaschke product}. Recall that a Blaschke product
 is called indestructible (see \cite{McL1972,Ros2008}) if 
for any conformal automorphism $T$ of the unit disk
the composition $T \circ B_{\mathcal{C}}$ is again a Blaschke product.
Note that postcomposition by such a conformal automorphism does not change the critical
set of a function in $H^{\infty}$.
Therefore, for any conformal automorphism $T$
of $\D$, we call the Blaschke product $T \circ B_\mathcal{{C}}$ a \textit{maximal
  Blaschke product with critical set $\mathcal{C}$}.

\medskip

If $N=0$, then the extremal problem (\ref{eq:ex1}) is exactly the problem of
maximizing the derivative at a point, i.e., exactly the character of Schwarz' lemma.
Let us put this observation in perspective.

\begin{remark}[The Nehari--Schwarz Lemma] \label{rem:1}
In the special case where $\mathcal{C}$ is a \textit{finite} sequence,
 Theorem \ref{thm:2} is essentially the classical and well--known Nehari--Schwarz lemma. 
\begin{itemize}
\item[(a)] 
In fact, if $\mathcal{C}=\emptyset$, then
  $\mathcal{F}_{\mathcal{C}}=H^{\infty}$, so all bounded analytic functions
  are competing functions, and
 Theoren \ref{thm:2}  is just the statement of Schwarz' lemma, which implies
 that $B_{\emptyset}$ is the identity map.
  In particular, the maximal Blaschke products without critical points, i.e.,
  the locally univalent maximal Blaschke products are precisely the
  unit disk automorphisms.
\item[(b)]
If $\mathcal{C}\not=\emptyset$ is a finite sequence and $N=0$, then Theorem \ref{thm:2} is exactly Nehari's
1947 generalization of Schwarz' lemma (see \cite{Neh1946}, 
Corollary\footnote{In his formulation of this Corollary, Nehari apparently
  assumes, implicitly, that the origin is not a critical point. Otherwise, Nehari's
statement concerning the case of equality would not be entirely correct.} to Theorem 1). In particular, if $\mathcal{C}=(z_1,\ldots, z_m)$ is
a finite sequence consisting of $m$ points, then every maximal Blaschke product with critical set $\mathcal{C}$
 is a finite Blaschke product of degree $m+1$. As we shall see in Remark \ref{rem:maxi1}
below, the converse is also true.
 Hence the maximal Blaschke products with finitely many critical points 
 are precisely the finite Blaschke products.
\end{itemize}
By these remarks, Theorem \ref{thm:2} might be viewed as  an extension of the
Nehari--Schwarz Lemma to arbitrarily many critical points.
\end{remark}

We now return to the extremal problem (\ref{eq:ex1}) and to a discussion of maximal Blaschke products and their
properties. As we shall see, maximal Blaschke products do have similar
characteristics as Bergman space inner functions on the one hand and display
many properties of finite
Blaschke products on the other hand. 

\medskip

We begin by relating maximal Blaschke products with Bergman space inner
functions. For this we note that a similar type of
extremal problem  as (\ref{eq:ex1}) 
was considered before for various
classes of analytic functions such as Hardy spaces and Bergman spaces, but
{\it with prescribed zeros instead of prescribed critical points}.
The following remark decribes this connection in full detail.

\begin{remark}[Hardy spaces and Bergman spaces]
If the sequence ${\cal C}$ is the critical set of a bounded analytic
function and if $N$ denotes the multiplicity of the point $0$ in $\mathcal{C}$, then 
according to Theorem \ref{thm:2}
 the maximal Blaschke product $B_{\mathcal{C}}$ with critical set ${\cal C}$
normalized by $B_{\mathcal{C}}(0)=0$ and $B^{(N+1)}_{\mathcal{C}}(0)>0$ is the unique solution to the
extremal problem


\begin{center}
\fbox{\parbox[h][\height -11mm +\baselineskip][c]{13cm}{
\begin{equation*}
 \max \big\{\Re f^{(N+1)}(0): f \in H^{\infty}\,, \, 
||f||_{\infty} \le 1 \text{ and } f'(z)=0 \text{ for } z \in{\cal C} \big \}
\, .
\end{equation*}
}}
\end{center}

\smallskip

This extremal property of a maximal Blaschke product  is reminiscent 
of a well--known extremal property of 
\begin{itemize}
\item[(i)]
Blaschke products in the Hardy spaces $H^{\infty}$ and
$$H^p:=\left\{ f : \D \to \C \text{ analytic} \,: {||f||}_p:=\left( \lim
    \limits_{r \to 1} 
\frac{1}{2\pi} \, \int
  \limits_{0}^{2\pi}  |f(r e^{it})|^p\,
dt  \right)^{1/p}<+ \infty \right\} \, , $$
where $1 <p  < +\infty$, 
\end{itemize}
as well as the
\begin{itemize}
\item[(ii)] canonical divisors in the
(weighted) Bergman spaces  
$$ {\cal A}_{\alpha}^p= \left\{
                 f:\D \to \C \text{ analytic} \,: {||f||}_{p,\alpha}:=\left(\frac{1}{\pi} \, \iint \limits_{\D} (1-|z|^2)^{\alpha} \, |f(z)|^p\,
d\sigma_z\right)^{1/p} <+ \infty \right\} \, , $$
where $-1<\alpha<+\infty$, $1 < p<+\infty$ and $d \sigma_z$ denotes
two--dimensional Lebesgue measure with respect to $z$,
\end{itemize} 
when the {\it zeros} are prescribed.

\medskip
More precsiely, let $1<p \le +\infty$, let
the sequence $ {\cal C}=(z_j)$ in $\D$ be the {\it zero set} of
a function in $H^p$ and let $N$ denote the multiplicity
of the point $0$ in ${\cal C}$.  Then the (unique) solution to the extremal problem


\begin{center}
\fbox{\parbox[h][\height -11mm +\baselineskip][c]{12cm}{
\begin{equation*}
\max\big\{\Re f^{(N)}(0): f \in H^p,\,  ||f||_p \le 1 \text{ and } f(z)=0
\text{ for } z \in{\cal C} \big\}
\end{equation*}
}}
\end{center}

\smallskip

is a Blaschke product $B$  with zero set ${\cal C}$ normalized by
$B^{(N)}(0)>0$, see \cite[\S 5.1]{DS}.

\medskip

 Hedenmalm \cite{Hed1991} (see also \cite{DKSS1993,DKSS1994}) had the idea of posing an
 appropriate counterpart of the latter extremal problem for 
Bergman spaces. His goal was to find a faithful analogue of Blaschke products
in Bergman spaces. 
As before, let ${\cal C}=(z_j)$
be a sequence in $\D$ where the point $0$ occurs $N$ times and assume that
$\mathcal{C}$ is the zero set of a function in $\mathcal{A}^p_{\alpha}$.
Then the extremal problem 


\begin{center}
\fbox{\parbox[h][\height -11mm +\baselineskip][c]{12cm}{
 \begin{equation*} 
  \max \big\{\Re f^{(N)}(0): f \in {\cal
  A}_{\alpha}^p\,, \,  || f||_{p,\alpha} \le 1 \text{ and } f(z)=0 \text{ for } z \in{\cal C}
\big \}\, 
\end{equation*}
}}
\end{center}

\smallskip

has a unique  extremal function ${\cal G} \in \mathcal{A}^p_{\alpha}$,  which vanishes
precisely on ${\cal C}$ and is normalized by ${\cal G}^{(N)}(0)>0$. 
The function $\mathcal{G}$ is called the canonical divisor for ${\cal C}$ or a
Bergman space inner function.
These functions play an extraordinary r$\hat{\text{o}}$le in the modern theory of Bergman
spaces, see \cite{HKZ,DS}.

\medskip

In summary, we have the following situation:

\smallskip

\begin{center}
\setlength{\extrarowheight}{6pt}
\begin{supertabular}{|l|c|l|}
\hline
{\bf prescribed data} & {\bf function space} & {\bf extremal function} \\[2mm]\hline
 critical set ${\cal C}$  & $H^{\infty}$ & maximal Blaschke product with
 critical set ${\cal C}$\\[2mm]\hline
 zero set ${\cal C}$       & $H^p$ & Blaschke product with zero set $\mathcal{C}$\\[2mm]\hline
zero set  ${\cal C}$       & ${\cal A}^p_{\alpha}$ & canonical divisor with
zero set
${\cal C}$ \\[2mm]
\hline
\end{supertabular}
\end{center}
\end{remark}

In light of this strong analogy, 
one expects that maximal Blaschke products enjoy similar properties as 
Blaschke products in $H^p$ spaces and canonical divisors in Bergman spaces,
with the critical points playing the r$\hat{\mbox{o}}$le of the zeros. 
An example is  analytic continuability. It is
a familiar result in $H^p$ theory that a Blaschke product 
has a holomorphic extension across every open arc of $\partial \D$ which does not
contain any limit point of its zero set, 
see \cite[Chapter II, Theorem 6.1]{Gar2007}. The same is true for a canonical
divisor in Bergman spaces. This
 was proved by
 Sundberg \cite{Sun1997} in 1997, who  improved earlier work
of Duren, Khavinson, Shapiro and Sundberg \cite{DKSS1993,DKSS1994} and Duren,
Khavinson and Shapiro \cite{DKS1996}.
Now keeping in mind that the critical points of maximal Blaschke products
take the r$\hat{\mbox{o}}$le of the zeros of  Blaschke products and  canonical divisors respectively, one
expects that a maximal Blaschke product has an analytic continuation across every open arc of
$\partial \D$ which does not meet any limit point of its critical set. This in fact
turns out to be true:

\begin{theorem}[Analytic continuability of maximal Blaschke products]\label{thm:analytic_continuation}
Let $F: \D \to \D$ be a maximal Blaschke product  with critical set ${\cal C}$. Then $F$ has an
analytic continuation across any arc of $\partial \D$ which is free of limit
points of ${\cal C}$.
\end{theorem}

Since a Blaschke product has an analytic continuation to a point $\xi
\in \partial \D$ if and only if $\xi$ is not a limit point of its zeros,
Theorem \ref{thm:analytic_continuation} leads to the following conclusion.

\begin{corollary}
The limit points of the critical set of a maximal Blaschke product coincide with the
limit points of its zero set.
\end{corollary}

We next turn to  a result about the structural properties of $H^{\infty}$
critical sets. It  follows from the results in \cite{Kra2011a}
 that the union of two $H^{\infty}$ critical sets is not necessarily
an $H^{\infty}$ critical set. However, if the two $H^{\infty}$ critical sets
have no common accumulation point on the unit circle, then their union is
again an $H^{\infty}$ critical set as the following result shows.

\begin{theorem} \label{thm:union}
Let $\mathcal{C}_1$ and $\mathcal{C}_2$ be two $H^{\infty}$ critical sets such
that $\overline{\mathcal{C}_1} \cap \overline{\mathcal{C}_2} \cap \partial
\D=\emptyset$. Then $\mathcal{C}_1 \cup \mathcal{C}_2$ is an $H^{\infty}$
critical set.
\end{theorem}

The analogous statement for the zero sets of Bergman space inner functions is
due to Sundberg \cite{Sun1997}.

\medskip

We next shift attention to maximal Blaschke products as generalizations of 
finite Blaschke products. Recall that in view of Remark \ref{rem:1} (b) every finite Blaschke product is a maximal
Blaschke product.  A rather strong property of the class of finite Blaschke products is their semigroup
property with respect to composition. In contrast, the composition of two
  Blaschke products does not need to be 
a Blaschke product (just consider non--indestructible Blaschke products). However, in case of
maximal Blaschke products the following result holds.

\begin{theorem}[Semigroup property of maximal Blaschke products] \label{thm:semigroup}
The set of maximal Blaschke products is closed under composition.
\end{theorem}

Finally, we consider the boundary behaviour of maximal Blaschke
products. Heins \cite{Hei86} showed that a function $B \in H^{\infty}$ is a
\textit{finite} Blaschke product if and only if
\begin{equation} \label{eq:bdd}
\lim \limits_{z \to \zeta} \left(1-|z|^2 \right)
\frac{|B'(z)|}{1-|B(z)|^2}=1
\end{equation}
for every $\zeta \in \partial \D$. The next theorem gives a partial extension
of this result for maximal Blaschke products.

\begin{theorem}[Boundary behaviour of maximal Blaschke products] \label{thm:bdd}
Let $B$ be a maximal Blaschke product with critical set $\mathcal{C}$. Then
(\ref{eq:bdd}) holds 
for every $\zeta \in \partial \D \backslash \overline{\mathcal{C}}$.
\end{theorem}

The results of the present paper are obtained by using the
 equivalence of two types of sets, critical sets for $H^{\infty}$ and zero sets for conformal
Riemannian pseudometrics  with curvature at most $-4$, see Corollary
\ref{cor:main0} below.  It turns out that the
latter zero sets are simpler to work with. Accordingly, we start in Section 2
by describing the relation between bounded analytic functions and negatively
curved conformal pseudometrics.  Section 3 contains the proofs of the main
results of this paper. It also gives a characterization of 
maximal Blaschke products in terms of ``maximal'' conformal pseudometrics.
In Section 4, we  generalize our results to analytic functions
defined on simply connected
proper subdomains of the complex plane other than the unit disk by 
using the Riemann mapping theorem. There, we  also indicate a connection
between maximal Blaschke products and the well--known 
Ahlfors' map for domains of finite connectivity $n \ge 2$. We close the paper
with a final  Section 5, which presents a number of open problems.

\bigskip

\textit{Acknowledgement.} We wish to thank an anonymous referee for carefully
reading the paper and providing us with a number of suggestions.

\section{Auxiliary results}

The proofs in this paper are based on  conformal
(Riemannian) pseudometrics and rely in particular on the results of \cite{Kra2011a}.
We first give a quick account of the relevant facts from conformal geometry and
refer to \cite{BM2007, Hei62,  KL2007, Krantz, Kra2011a, KR2008, Smi1986}
for more information.

\medskip

We call a continuous nonnegative function $\lambda : G \to \R$ defined on a
domain $G \subseteq \C$ a conformal density   and the quantity $\lambda(z) \,
|dz|$ a conformal pseudometric on $G$.
 We say
$\lambda(z)\, |dz|$ has a zero of order $m_0>0$ at $z_0 \in G$ if
\begin{equation*}
\lim_{z \to z_0} \frac{\lambda(z)}{|z-z_0|^{m_0}} \quad\text{ exists and }
\not=0 \, . 
\end{equation*}
In this paper we will only consider conformal
pseudometrics $\lambda(z) \, |dz|$  with isolated
zeros.
 A  sequence ${\cal C}=(\xi_j) \subset G$
\begin{equation*}
(\xi_j):=(\underbrace{z_1, \ldots, z_1}_{m_1  -\text{times}},\underbrace{z_2,
  \ldots, z_2}_{m_2-\text{times}} , \ldots  )\,, \, \,  z_k \not=z_n \text{ if
} k\not=n, 
\end{equation*}
is called  the zero set of a conformal
pseudometric $\lambda(z) \, |dz|$, if $\lambda(z)>0$ for $z \in
G\backslash {\cal C}$ and if $\lambda(z)\, |dz|$ has a zero of order $m_k \in \N$
at $z_k$ for all $k$. 
 We will always assume that  $\lambda$ is of class
$C^2$ in a neighborhood of any point $z_0 \in G$ where $\lambda(z_0)>0$.
Hence the curvature $\kappa_{\lambda}$
of $\lambda(z) \, |dz|$ can be defined by
\begin{equation} \label{def:curvature}
 \kappa_{\lambda}(z_0)=-\frac{\Delta (\log \lambda)}{\lambda^2}(z_0) \, 
\end{equation}
 for any point $z_0\in G$ with $\lambda(z_0)>0$.

\medskip

An important aspect of  curvature is its conformal invariance. Suppose $f:G \to D$ is an analytic function from a
domain $G$ into a domain $D$ and let $D$ be equipped with a positive conformal
pseudometric $\lambda(z)\, |dz|$ with curvature $\kappa_{\lambda}$. Then the
pullback to $G$ via $f$ of $\lambda(z)\, |dz|$ is a conformal pseudometric on
$D$ defined by
\begin{equation*}
f^*\lambda(z)\, |dz|:= \lambda(f(z))\, |dz|\,.
\end{equation*}
Now  the curvature of the pullback pseudometric $f^*\lambda(z)\, |dz|$ and the
pseudometric $\lambda(w) \, |dw|$ are related by the fundamental identity
\begin{equation*}
\kappa_{f^*\lambda}(z)=\kappa_{\lambda}(f(z)) \, ,
\end{equation*}
which is valid for every $z \in G \backslash\{z\in G\, :\,  f'(z)=0\}$. Note that in particular
the zero set
of the pseudometric $f^*\lambda(z)\, |dz|$ is precisely the critical set of
the function $f$.

\medskip

The ubiquitous example for a conformal  pseudometric is the Poincar\'e metric or
hyperbolic metric
\begin{equation*}
\lambda_{\D}(z)\, |dz| =\frac{1}{1-|z|^2}\, |dz|
\end{equation*}
for the unit disk $\D$; it has constant curvature $-4$ on $\D$. The hyperbolic
metric $\lambda_{\D}(z) \, |dz|$ has the following
important property.

\begin{theorem}\label{thm:ahlfors}
Let $\lambda(z)\,|dz|$ be a conformal pseudometric on
$\D$ with curvature $\kappa_{\lambda}(z) \le -4$ for all $z\in \D$ with
$\lambda(z)>0$. Then the following statements hold.
\begin{itemize}
\item[(a)] For every $z \in \D$
 $$\lambda(z) \le \lambda_{\D}(z)\, .$$
\item[(b)]
If $\lambda(z_0)=\lambda_{\D}(z_0)$ for some $z_0\in \D$, then
$\lambda(z)=\lambda_{\D}(z)$ for all $z \in \D$.
\end{itemize}
\end{theorem}

Theorem \ref{thm:ahlfors} (a) is due to Ahlfors \cite{Ahl1938} and it is usually
called Ahlfors' lemma. The case of equality in Ahlfors' lemma, i.\;\!e.~Theorem
\ref{thm:ahlfors} (b), was proved by
Heins \cite[\S7]{Hei62}, see also Royden \cite{Roy1986} and Minda \cite{Min1987}.

\medskip

The following theorem gives a sharpening of Theorem
\ref{thm:ahlfors} for  conformal pseudometrics with prescribed zero set.

\begin{theorem}[Maximal conformal pseudometric] \label{thm:hilf1}
Let $\mathcal{C}$ be a sequence of points in $\D$ such that there exists
 a conformal pseudometric in $\D$ with zero set $\mathcal{C}$ and
curvature $\le -4$ in $\D \backslash
\mathcal{C}$.  Then there exists a unique
conformal pseudometric $\lambda_{max}(z) \, |dz|$ on $\D$ with zero set $\mathcal{C}$
and constant curvature $-4$ on $\D \backslash \mathcal{C}$ such that
for any conformal pseudometric $\lambda^*(z) \, |dz|$ with zero set $\mathcal{C}^*
\supseteq \mathcal{C}$ and curvature $\kappa_{\lambda^*}(z) \le -4$ on $\D
\backslash \mathcal{C}^*$ the following conditions hold.
\begin{itemize}
\item[(a)] For every $z \in \D$
$$ \lambda^*(z) \le \lambda_{max}(z) \, .$$
\item[(b)] If 
$$ \lim \limits_{z \to z_0} \frac{\lambda^*(z)}{\lambda_{max}(z)} =1$$
for some point $z_0 \in \D$, then $\lambda^*(z)=\lambda_{max}(z)$ for all $z
\in \D$.
\end{itemize}
\end{theorem}

We note that the first statement of Theoem \ref{thm:hilf1} is a result by Heins
\cite[Theorem 13.1]{Hei62}. It suggests calling the conformal pseudometric
$\lambda_{max}(z)\, |dz|$ the maximal conformal pseudometric on $\D$ with zero
set ${\cal C}$ and curvature $-4$ on $\D\backslash {\cal C}$.

\medskip

Theorem \ref{thm:hilf1} (b) will play a crucial r$\hat{\mbox{o}}$le in the proof
of Theorem \ref{thm:2}. In order to prove it, we need the following result.

\begin{lemma}\label{thm:existence1}
Let $G$ be a bounded and regular domain\footnote{i.\!\;e.~there exists a
  Green's function for $G$ which vanishes
  continuously on the boundary of $G$.}, let $b: \partial G \to (0, +
\infty)$ be a continuous function on the boundary $\partial G$ of $G$
 and let $\kappa$ be a bounded, nonpositive and locally
H\"older continuous on $G$. Then there exists a unique positive conformal
pseudometric $\lambda(z)\, |dz|$ on $G$ with curvature $\kappa_{\lambda}\equiv
\kappa$ on $G$ such that $\lambda$ is continuous on
the closure
$\overline{G}$ and
$\lambda(\xi)=b(\xi)$ for all $\xi \in \partial G$.
\end{lemma}

For completeness, we sketch a proof of Lemma \ref{thm:existence1} here using a
standard result from nonlinear elliptic PDEs.

\medskip

{\bf Proof.} We first note that the Dirichlet problem
\begin{equation}
\begin{array}{ll} \label{eq:di}
\Delta u=-\kappa(z) \, e^{2u} & \quad \text{ if } z \in G \, , \\
u(\xi) =\log b(\xi) & \quad \text{ if } \xi \in \partial G \, .
\end{array}
\end{equation}
has a unique solution $u \in C^2(G) \cap C(\overline{G})$, see 
 \cite[p.~53--55 \& p.~304]{GT} and \cite[p.~286]{Cou1968}.
Taking into account the definition of curvature, i.e.~the formula (\ref{def:curvature}), 
this means that
$$ \lambda(z) \, |dz|:=e^{u(z)} \, |dz|$$
is  the unique positive  conformal pseudometric on $G$ with curvature
$\kappa$ such that $\lambda$ is continuous on $G$ and $\lambda(\xi)=b(\xi)$ for every
$\xi \in \partial G$.
\hfill{$\blacksquare$}

\begin{remark} \label{rem:max}
The uniqueness statement in Lemma \ref{thm:existence1} ultimately comes from
the maximum principle for solutions to the Dirichlet problem (\ref{eq:di}),
see \cite[Theorem 10.1]{GT}. In the terminology of conformal pseudometrics, the maximum principle says
the following.
Let $\lambda(z) \, |dz|$ and $\mu(z) \, |dz|$ be positive conformal
pseudometrics on a bounded domain $G$ with the following properties:
$\kappa_{\lambda} \le \kappa_{\mu} \le 0$
 on $G$, $\lambda$ and $\mu$ are continuous and positive on $\overline{G}$ and
$\lambda(\xi) \le \mu(\xi)$ for all $ \xi \in \partial G$. Then $\lambda \le
\mu$ throughout $G$. 
\end{remark}

Before passing to the proof of Theorem \ref{thm:hilf1} (b), we note that the
special case of Theorem \ref{thm:hilf1} (b) as described in Theorem
\ref{thm:ahlfors} (b) allows for a quick proof using the so--called strong
maximum principle of E.~Hopf 
from the theory of elliptic PDEs, see e.g.~\cite{GT,Min1987}. In the general
case, we have to deal with the problem that the assumption
$$ \lim \limits_{z \to z_0} \frac{\lambda^*(z)}{\lambda_{max}(z)} =1$$
in Theorem \ref{thm:hilf1} (b) might hold for a point $z_0 \in \D$ where
$\lambda_{max}$ vanishes.  As we shall see, Lemma \ref{thm:existence1}
will enable us to still make use of Hopf's strong maximum principle.

\medskip

{\bf Proof of Theorem \ref{thm:hilf1} (b).}
We divide the proof in two cases.

\smallskip

(i) \, First assume that $\lambda^*(z_1)=\lambda_{max}(z_1)$ for some point
$z_1 \in G:=\D \backslash \mathcal{C}^*$.
Then the function 
$$ u(z):=\log \frac{\lambda_{max}(z)}{\lambda^*(z)}$$
is twice continuously differentiable and nonnegative on $G$ and
\begin{eqnarray*} 
 \Delta u (z) &=&\Delta \log \lambda_{max}(z)-\Delta \log \lambda^*(z) \le 4
\lambda_{max}(z)^2-4 \lambda^*(z)^2\\ &=& 4 \lambda_{max}(z)^2 (1-e^{-2
  u(z)})\le 8 \lambda_{max}(z)^2 \, u(z) \, 
\end{eqnarray*}
for every $z \in G$.
The strong maximum principle (\cite[Theorem 3.5]{GT}) implies that $u \equiv 0$ in $G$,
i.e.~$\lambda^*(z)=\lambda_{max}(z)$ for all $z \in \D$ by continuity.

\smallskip

(ii) \, Now consider the case that $\lambda^*(z)<\lambda_{max}(z)$ for all $z \in
\D \backslash \mathcal{C}^*$. This implies $\lambda^*(z)<\lambda_{max}(z)$ for all $z \in
\D \backslash \mathcal{C}$. In particular, the point
 $z_0$ belongs to $\mathcal{C}$. Let $m$ denote the order of the zero
 of $\lambda_{max}(z) \, |dz|$ at $z_0$.  Since $\mathcal{C}$ is a
 discrete subset of the unit disk $\D$,  there exists an open disk
 $K$ compactly contained in $\D$ such that $z_0 \in K$,  $\lambda^*(z)>0$ for
 every $z \in\overline{K} \backslash \{z_0\}$ and
\begin{equation} \label{eq:ha}
\lambda^*(\xi)<\lambda_{max}(\xi)  \quad \text{ for
 all } \xi \in \partial K \, .
\end{equation}
 We consider the pseudometric
$$ \tilde{\lambda}_{max}(z)\, |dz|:=\frac{\lambda_{max}(z)}{|z-z_0|^m} \,
|dz| \, $$
on the punctured disk $K \backslash \{ z_0\}$. Clearly, $\tilde{\lambda}_{max}$ is
of class $C^2$ in $K \backslash \{ z_0\}$ and a
 quick computation shows that $\tilde{\lambda}_{max}(z)\, |dz|$  
has curvature $\kappa(z)=-4 |z-z_0|^{2m}$ there. Since $z_0$ is a zero of
$\lambda_{max}(z) \, |dz|$ of order $m$, the density  
 $\tilde{\lambda}_{max}$  has a continuous extension to $K$.
We claim that (the continuous extension of) $\tilde{\lambda}_{max}$ is of class $C^2$ on the entire disk $K$.
In fact, Lemma \ref{thm:existence1} shows that there is a positive conformal
pseudometric $\mu(z)\, |dz|$  on $K$ with curvature $\kappa$ there,
continuous on $\overline{K}$  and
$\mu(\xi)=\tilde{\lambda}_{max}(\xi)$ for every $\xi \in\partial K$.
In particular, $\mu$ is of class $C^2$ on $K$. It now turns out that $\mu \equiv \tilde{\lambda}$.
 To
see this, consider the auxiliary function
$$ v(z):=\max \left\{ 0, \log \frac{\tilde{\lambda}_{max}(z)}{\mu(z)} \right\} \, , \qquad z
\in K \, .$$
If $z_1 \in K$ such that $v(z_1)>0$, then 
$$\Delta v(z_1)=\Delta \log \tilde{\lambda}_{max}(z_1)-\Delta \log \mu(z_1)
=4 |z_1-z_0|^{2m} \tilde{\lambda}_{max}(z_1)^2-4 |z_1-z_0|^{2m} \mu(z_1)^2
 \ge 0\, , $$ so $v$ is a continuous subharmonic
function on $\{z \in K \, : \, v(z)>0\}$ with boundary values $0$, i.e., $v
\le 0$ by the maximum principle for subharmonic functions. This shows that $v \equiv 0$, i.e.,
~$\tilde{\lambda}_{max} \le \mu$. A similar argument establishes the reverse
inequality, so we have proved that $\tilde{\lambda}_{max}=\mu$ is of class
$C^2$ on the entire disk $K$.

\smallskip

We now go back to inequality (\ref{eq:ha}).
Lemma \ref{thm:existence1} shows that there exists
 a positive pseudometric $\nu (z) \, |dz|$ on $K$ with curvature $\kappa(z)=-4 |z-z_0|^{2m}$ on $K$
 such that $\nu$ is continuous on $\overline{K}$ and 
$$ \frac{\lambda^*(\xi)}{|\xi-z_0|^m}
<\nu(\xi)<\frac{\lambda_{max}(\xi)}{|\xi-z_0|^m} =\tilde{\lambda}_{max}(\xi) \quad \text{ for every } \xi
\in \partial K \, .$$
Then 
\begin{equation} \label{eq:hu}
\frac{\lambda^*(z)}{|z-z_0|^m}\le \nu(z) \le \tilde{\lambda}_{max}(z) \quad \text{ for
all } z \in K \backslash \{ z_0\}\, .
\end{equation}
The right inequality follows directly from the maximum principle (see Remark \ref{rem:max})
applied to the conformal pseudometrics $\nu(z) \, |dz|$ and
$\tilde{\lambda}_{max}(z) \, |dz|$. In order to prove the left inquality in
(\ref{eq:hu}), consider the function $s : K \backslash \{ z_0\} \to \R$
defined by
$$ s(z):=\max \left\{ 0, \log \left(\frac{\lambda^*(z)}{|z-z_0|^m \nu(z)} \right) \right\} \,  .$$
If $z \in K\backslash \{z_0\}$ such that $s(z)>0$, then
$$ \Delta s(z)=\Delta \log \lambda^*(z)-\Delta \log \nu(z) \ge 4
\lambda^*(z)^2-4 |z-z_0|^{2m} \nu(z)^2 \ge 0 \, .$$
Hence, the nonnegative function $s : K \backslash \{ z_0\} \to \R$ is subharmonic on $K
\backslash \{ z_0\}$. Moreover, since $\lambda^*(z)$ has a zero at $z_0$ of at
least order $m$, the function $s$ has a continuous extension to the point
$z_0$. As a consequence, the function $s$ has a subharmonic extension to the entire disk
$K$ with boundary values $0$. Thus the maximum principle for subharmonic functions
leads to $s \le 0$ throughout $K$, so, by the definition of the function $s$, the left inquality in
(\ref{eq:hu}) follows.

\smallskip

We next combine the assumption of Theorem \ref{thm:hilf1} (b) with the
estimates (\ref{eq:hu}) and obtain
$$ \lim \limits_{z \to z_0} \frac{\tilde{\lambda}_{max}(z)}{\nu(z)}=\lim \limits_{z \to z_0} \frac{\lambda_{max}(z)}{\lambda^*(z)}=
1 \, .$$
Since, as we have observed above, the functions $\nu$ and $\tilde{\lambda}_{max}$ are of class $C^2$ on $K$,
this means that the auxiliary function
$$\tilde{u}(z):=\log  \frac{\tilde{\lambda}_{max}(z)}{\nu(z)}$$
is nonnegative and of class $C^2$ in $K$ with $\tilde{u}(z_0)=0$. Now,
similar to part (a), one can show that
$$ \Delta \tilde{u}(z)\le 8 \lambda_{max}(z)^2 |z-z_0|^{2m} \, \tilde{u}(z) \, , \qquad
z \in K \, . $$
Hence the strong maximum principle gives $\tilde{u} \equiv 0$, contradicting the boundary
condition $\tilde{u}(\xi)>0$ for every $\xi \in \partial K$.
We conclude that our assumption $\lambda^*\not\equiv \lambda_{max}$ cannot
hold and the proof of Theorem \ref{thm:hilf1} (b) is complete.
\hfill{$\blacksquare$}

\medskip

The next result shows that we can represent the maximal pseudometric
$\lambda_{max}(z) \, |dz|$ in
Theorem \ref{thm:hilf1} as the pullback of the Poincar\'e metric
$\lambda_{\D}(z) \, |dz|$ under a
specific function $F \in H^{\infty}$.

\begin{theorem} \label{thm:hilf2}
Let $\mathcal{C}$ be an $H^{\infty}$ critical set and let $\lambda_{max}(z) \,
|dz|$ be the maximal conformal pseudometric on $\D$ with zero
set $\mathcal{C}$ and curvature $-4$ on $\D \backslash \mathcal{C}$.
 Then the following statements hold.
\begin{itemize}
\item[(a)] There exists a function $F \in H^{\infty}$ with critical set
  $\mathcal{C}$ such that
\begin{equation} \label{eq:pullmax}
 \lambda_{max}(z)=\frac{|F'(z)|}{1-|F(z)|^2} \, , \qquad z \in \D \, .
\end{equation}
Moreover, $F$ is uniquely determined by $\lambda_{max}(z) \, |dz|$ up to postcomposition
with a unit disk automorphism.
\item[(b)] The function $F$ in (a) is an indestructible Blaschke product with
  critical set $\mathcal{C}$.
\end{itemize}
\end{theorem}

\begin{remark}\label{rem:wu}
\begin{itemize}
\item[(a)] Statement (a) and part of statement (b) 
 of Theorem \ref{thm:hilf2} can be found in the work of Heins
\cite[Theorem 29.1]{Hei62}. His proof is in three steps. In a first step,
Heins  proved  that for a finite sequence
$\mathcal{C}$ (consisting of $m$ points) in $\D$ there is always a finite Blaschke
product of degree $m+1$ with critical set $\mathcal{C}$ and that this finite Blaschke
product is uniquely determined by its critical set $\mathcal{C}$ up to postcomposition by a unit
disk automorphism. His second step 
consists of showing that if $B$ is a finite Blaschke product (of degree $m+1$)
with critical set $\mathcal{C}$ (then containing exactly $m$ points,
cf.~\cite[p.~78]{SS}) then the
pullback  of $\lambda_{\D}(z) \, |dz|$ via $B$,
$$ \frac{|B'(z)|}{1-|B(z)|^2} \, |dz|$$
is the maximal conformal pseudometric $\lambda_{max}(z) \, |dz|$ on $\D$ with zero set $\mathcal{C}$ and
curvature $-4$ on $\D \backslash \mathcal{C}$. Finally, part (a) of Theorem
\ref{thm:hilf2} is established by letting $m$ tend to $\infty$.

\medskip

We note that the first and the second step together prove Theorem \ref{thm:hilf2} (b)
for {\it finite} sequences $\mathcal{C}$.
The general case of Theorem \ref{thm:hilf2} (b) is more involved and is proved in \cite{Kra2011a}.
\item[(b)] 
 Theorem \ref{thm:hilf2} (a) is  actually a special instance of ``Liouville's Theorem'':
\begin{quote} {\it  
Every conformal pseudometric $\lambda(z) \, |dz|$ in
$\D$ with zero set $\mathcal{C}$ and constant curvature $-4$ on $\D \backslash \mathcal{C}$ can be written as
the pullback of the hyperbolic metric $\lambda_{\D}(w) \, |dw|$ under a function $f\in H^{\infty}$ 
with critical set $\mathcal{C}$, that is
\begin{equation} \label{eq:pullback}
 \lambda(z)=\frac{|f'(z)|}{1-|f(z)|^2} \, , \qquad z \in \D \, .
\end{equation}
 Moreover, the analytic function $f : \D
\to \D$ (the so--called developing map of $\lambda(z) \, |dz|$) is uniquely
determined by the pseudometric $\lambda(z) \, |dz|$ up to postcompositon
with a unit disk automorphism. } \end{quote}

In fact,  Liouville  \cite{Lio1853} stated only the
zerofree case ($\mathcal{C}=\emptyset$) and his proof is not entirely convincing by
today's standards. 
A very elegant geometric proof of the zerofree case was given by D.~Minda
in his notes \cite{Min}.
 The general case of Liouville's Theorem ($\mathcal{C}\not=\emptyset$)
can be found e.g.~in \cite{Bie16, CW94, CW95, KR2008, Nit57, Yam1988}.  
Note that in Liouville's theorem, the zeros of the pseudometric $\lambda(z) \,
|dz|$ are precisely the critical points of its developing map $f \in
H^{\infty}$.
\item[(c)] Using the terminology of Liouville's theorem,
Theorem \ref{thm:hilf2} (b)  says that for a sequence $\mathcal{C}$ in $\D$
the developing maps of the maximal conformal pseudometric in $\D$
with zero set $\mathcal{C}$ and curvature $-4$ in $\D \backslash \mathcal{C}$
are indestructible Blaschke products with critical set $\mathcal{C}$.
\end{itemize}
\end{remark}

The next result provides a simple criterion for a sequence $\mathcal{C}$ being
 an $H^{\infty}$ critical set, which will be useful later on.

\begin{corollary} \label{cor:main0}
Let $\mathcal{C}$ be a sequence of points in $\D$. Then the following
conditions are equivalent.
\begin{itemize}
\item[(a)] $\mathcal{C}$ is an $H^{\infty}$ critical set.
\item[(b)] $\mathcal{C}$ is the zero set of a conformal pseudometric in $\D$
  with zero set $\mathcal{C}$ and curvature $\le -4$ in 
$\D \backslash \mathcal{C}$.
\end{itemize}
\end{corollary}

{\bf Proof.}
Let $\mathcal{C}$ be the critical set of $f \in H^{\infty}$. Then the pullback
$$ f^*\lambda_{\D}(z) \, |dz|=\frac{|f'(z)|}{1-|f(z)|^2} \, |dz|$$
of the hyperbolic metric $\lambda_{\D}(w) \, |dw|$ under $f$ is a conformal
pseudometric with zero set $\mathcal{C}$ and curvature 
 $-4$ on $\D \backslash \mathcal{C}$.
This proves ``(a) $\Longrightarrow$ (b)''. Conversely, if (b) holds, then
Theorem \ref{thm:hilf1} implies that there is a maximal conformal pseudometric
$\lambda_{max}(z) \, |dz|$ with zero set $\mathcal{C}$ and curvature $-4$ on
$\D \backslash \mathcal{C}$. Using Theorem \ref{thm:hilf2} (a), we see that
there is a function $F \in H^{\infty}$ with critical set $\mathcal{C}$. 
\hfill{$\blacksquare$}

\begin{remark}
Heins \cite[\S25 \& \S 26]{Hei62} initiated the study of the mapping properties of the developing
maps  of maximal conformal pseudometrics. He obtained
  some necessary conditions as well as sufficient conditions   
for  developing maps of maximal conformal pseudometrics,
 but he did not prove that they are always
Blaschke products. He also posed the 
problem of characterizing  the developing maps of maximal conformal pseudometrics, cf.~\cite[\S 26 \& \S 29]{Hei62}.
\end{remark}

\section{Proofs}

\subsection{Proof of Theorem \ref{thm:2} and some consequences}

The proof of Theorem \ref{thm:2} consists in identifying the extremal function(s) for the extremal
problem (\ref{eq:ex1}) as the developing maps of the {\it maximal} conformal
pseudometric $\lambda_{max}(z) \, |dz|$ with zero set $\mathcal{C}$ and
curvature $-4$ on $\D \backslash \mathcal{C}$. This is accomplished with the help of
Theorem \ref{thm:hilf1} (b) and an application of
Theorem \ref{thm:hilf2} (b).

\bigskip

{\bf Proof of Theorem \ref{thm:2}.}

(a) \, We first note that a normal family argument ensures the existence of an
extremal function for the extremal problem (\ref{eq:ex1}), i.e., there is an
analytic function  $g \in \mathcal{F}_{\mathcal{C}}$ such that $\Re g^{(N+1)}(0) \ge \Re f^{(N+1)}(0)$
for all $f \in \mathcal{F}_{\mathcal{C}}$. It is easy to show that
$g^{(N+1)}(0)$ is real and positive and that $g(0)=0$.

\medskip

In fact, since $\mathcal{C}$ is an
$H^{\infty}$ critical set and $N$ is the number of times that $0$ occurs in the
sequence $\mathcal{C}$, there exists some $f \in \mathcal{F}_{\mathcal{C}}$
with $f^{(N+1)}(0)\not=0$. In view of
 $\eta f \in \mathcal{F}_{\mathcal{C}}$ for any $f \in
\mathcal{F}_{\mathcal{C}}$ and any $\eta \in \partial \D$, this implies
$g^{(N+1)}(0)=\Re g^{(N+1)}(0)>0$. 
Since
\begin{equation*}
 \tilde{g}(z):= \frac{g(z)-g(0)}{1-\overline{g(0)} \, g(z)}
= \frac{g^{(N+1)}(0)}{1-|g(0)|^2} \frac{1}{(N+1)!}\, z^{N+1}+\cdots 
   \end{equation*}
 belongs to ${\cal F}_{\cal C}$, we deduce that $g(0)=0$.

\medskip

(b) Let
$$ \lambda^*(z) \,|dz|:=\frac{|g'(z)|}{1-|g(z)|^2}\,|dz| \, , \qquad z\in \D \, ,$$
be the pullback of the hyperbolic metric $\lambda_{\D}(z) \, |dz|$ under the
holomorphic function $g$.
Then $\lambda^*(z) \, |dz|$ is a conformal pseudometric on $\D$ with zero set
  $\mathcal{C}^* \supseteq \mathcal{C}$ and curvature $-4$ on $\D \backslash
  \mathcal{C}^*$. By Theorem \ref{thm:hilf1} there exists a maximal conformal pseudometric  $\lambda_{max}(z) \, |dz|$
on $\D$ with zero set $\mathcal{C}$ and curvature $-4$ in $\D \backslash
\mathcal{C}$. Let
$F \in H^{\infty}$ denote  the developing map of $\lambda_{max}(z) \, |dz|$ 
normalized by $F(0)=0$ and $F^{(N+1)}(0) \ge 0$. Since $\mathcal{C}$ is the
critical set of $F$, it follows that
$F^{(N+1)}(0) > 0$ and $F \in {\cal F}_{\cal C}$, so 
\begin{equation} \label{eq:help1}
0<F^{(N+1)}(0)\le g^{(N+1)}(0) \, .
\end{equation}
 On the other hand, the maximal property of $\lambda_{max}(z) \, |dz|$  shows
\begin{equation*}
 \lambda^*(z) =\frac{|g'(z)|}{1-|g(z)|^2} \le
\frac{|F'(z)|}{1-|F(z)|^2}=\lambda_{max}(z)\, , \quad z \in \D\, ,
\end{equation*}
so
\begin{equation} \label{eq:help2}
\frac{g^{(N+1)}(0)}{F^{(N+1)}(0)}= \lim\limits_{z \to 0}
\frac{\lambda^*(z)}{\lambda_{max}(z)} \le 1\, .\end{equation}
Conditions (\ref{eq:help1}) and (\ref{eq:help2}) together imply
 $F^{(N+1)}(0)=g^{(N+1)}(0)$ and, appyling (\ref{eq:help2}) again, we get that
 $$\lim \limits_{z \to
  0} \frac{\lambda^*(z)}{\lambda_{max}(z)}=1 \, .$$ Hence  Theorem \ref{thm:hilf1} (b) 
gives $\lambda^* \equiv \lambda_{max}$ and Theorem \ref{thm:hilf2} (a)
 shows that $g=T \circ F$ for some conformal automorphism  $T$ of the unit disk. Since
 $g(0)=F(0)$ and $g^{(N+1)}(0)=F^{(N+1)}(0)$, we finally arrive at
 $g=F$. In particular, $g$ is uniquely determined.

\medskip 

(c) Theorem \ref{thm:hilf2} (b) shows that $g=F$ is an indestructible Blaschke product.
If $\mathcal{C}$ is a finite sequence consisting of $m$ points, then by a result of
Heins \cite{Hei62} (see Remark \ref{rem:wu} (a)), 
the function $F$ in part (b) above is  a finite
Blaschke product of degree $m+1$. If $\mathcal{C}$ is an infinite
sequence, then $F$ cannot be a finite Blaschke product, since finite Blaschke
products  only have finitely many critical points, see e.g.~\cite[p.~78]{SS}.
\hfill{$\blacksquare$}

\bigskip

We note the following immediate consequence of the above proof of Theorem \ref{thm:2}.

\begin{corollary} \label{cor:main1}
Let $\mathcal{C}$ be a sequence of points in $\D$  and $F \in H^{\infty}$ with
critical set $\mathcal{C}$. Then the following are equivalent.
\begin{itemize}
\item[(a)] $F$ is a maximal Blaschke product with critical set $\mathcal{C}$.
\item[(b)] The conformal pseudometric
$$ (F^*\lambda_{\D})(z) \, |dz|:=\frac{|F'(z)|}{1-|F(z)|^2} \, |dz|$$
is the maximal conformal pseudometric $\lambda_{max}(z) \, |dz|$ on $\D$ with  zero set $\mathcal{C}$ and
curvature $-4$ on $\D \backslash \mathcal{C}$.
\end{itemize}
\end{corollary}

\begin{remark} \label{rem:maxi1}
As mentioned earlier (see Remark \ref{rem:1} (a)),  any finite Blaschke product is a maximal Blaschke product. 
To see this, let  $F$ be a finite Blaschke product of degree $m+1$ say, then
$F$ has a finite critical set $\mathcal{C}=(z_1, \ldots, z_m)$.  
By Remark \ref{rem:wu} (a), $(F^*\lambda_{\D})(z) \, |dz|$ is the maximal conformal pseudometric on $\D$ with  zero set $\mathcal{C}$ and
curvature $-4$ on $\D \backslash \mathcal{C}$. So, Corollary \ref{cor:main1} shows
that $F$ is a maximal Blaschke product with critical set $\mathcal{C}$.
\end{remark}

The next result  follows easily from Theorem \ref{thm:2}. It  will
be very useful later.

\begin{corollary} \label{cor:main}
Let $(\mathcal{C}_n)$ be a sequence of $H^{\infty}$ critical sets with
$\ldots \mathcal{C}_{n+1} \supseteq \mathcal{C}_{n} \supseteq \ldots \supseteq
\mathcal{C}_1$ and let $F_n$ denote the extremal function for the extremal
problem (\ref{eq:ex1}) for $\mathcal{C}_n$. Assume that
$$ \mathcal{C}:=\bigcup \limits_{n=1}^{\infty} \mathcal{C}_n$$
is an $H^{\infty}$ critical set and let $F$ be the extremal function for the extremal
problem (\ref{eq:ex1}) for $\mathcal{C}$. Then the sequence $(F_n)$ converges to $F$ locally uniformly
in $\D$.
\end{corollary}

{\bf Proof.} Let $N_n$ denote the number of times that $0$ appears in the
sequence $\mathcal{C}_n$. From $\mathcal{C}_{n+1} \supseteq \mathcal{C}_n$, we
get $N_{n+1} \ge N_n$. Since $\mathcal{C}=\cup_{n=1}^{\infty} \mathcal{C}_n$ is an
$H^{\infty}$ critical set, we can  assume that $0$ occurs finitely often, say $N$ times in the sequence
$\mathcal{C}$. This implies that $N_n=N$ for all but finitely many $n$, say
for all $n \ge K$. 
Since $F \in \mathcal{F}_{\mathcal{C}_n}$ for every positive integer $n$, we deduce
 $\Re F^{(N+1)}(0) \le \Re F_{n}^{(N+1)}(0) \le \Re F_{n+1}^{(N+1)}(0)$ for all $n \ge K$.
If $g$ is a subsequential limit function of the sequence $(F_n)$  with respect to locally uniform convergence in
$\D$, then this implies $\Re g^{(N+1)}(0) \ge \Re F^{(N+1)}(0)$. On the other
hand, $g \in \mathcal{F}_{\mathcal{C}}$, so $\Re  g^{(N+1)}(0) \le \Re
F^{(N+1)}(0)$. By the uniqueness statement in Theorem \ref{thm:2}, we deduce $g=F$.
Since $(F_n)$ is a normal family, we conclude that $(F_n)$ has a unique
subsequential limit function. This proves the corollary.
\hfill{$\blacksquare$}

\subsection{Proof of Theorem \ref{thm:analytic_continuation}}

We now shift attention to the proof of Theorem
\ref{thm:analytic_continuation}. 

\bigskip

{\bf Proof of Theorem \ref{thm:analytic_continuation}.}
Let  $m_j$ denote the number of times that $z_j$ appears in the sequence ${\cal C}=(z_j)$. For simplicity we set
$z_1=0$ and
$N=m_1$. 
We may assume that the maximal Blaschke product $F$ for $\mathcal{C}$ is
normalized by $F(0)=0$ and  $F^{(N+1)}(0)>0$, so $F$ is the unique
extremal function for the extremal problem of Theorem \ref{thm:2}.
Suppose that $z_0 \in \partial \D$ is not a limit point of ${\cal C}$.
Then there exists an open disk $K$ with the following three properties:
$z_0 \in K$, the boundary $\partial K$ intersects
the unit circle perpendicularly and $\D \cap K$ contains none of the points of
${\cal C}$. We set $I:=K \cap \partial \D$. It suffices to
prove that $F$ has an analytic continuation to $\D \cup K$.

\medskip

Let $F_n$ denote the extremal function for the extremal problem (\ref{eq:ex1})
for the finite critical set
\begin{equation*}
{\cal C}_n:=\big(\underbrace{0,\ldots,0}_{N\text{--times}},\,
\underbrace{z_2,\ldots,z_2}_{m_2\text{--times}}, \ldots ,\, \underbrace{z_n,
  \ldots, z_n}_{m_n\text{--times}} \big) \, .\end{equation*} 
 Note that for any positive integer $n$ the function $F_n$ is
a finite Blaschke product of degree $N+m_2+ \cdots+ m_n +1$ 
with $F_n(0)=0$ and $F_n^{(N+1)}(0)>0$. Since by construction
$\cup_{n=1}^{\infty} \mathcal{C}_n=\mathcal{C}$, we see 
from Corollary \ref{cor:main}  that
$F_n \to F$ locally uniformly in $\D$.
We now consider the auxiliary functions
\begin{equation*}
 \varphi_n(z):=\frac{F_n(z)}{z F_n'(z)} 
 \end{equation*}
for $n=1,2, \ldots$.
Since each $F_n$ is a finite Blaschke product, the function
$\varphi_n$ is  analytic in $K$ and $0<\varphi_n(\zeta)  \le 1$ for each $\zeta
\in I$ by the Julia--Wolff--Carath\'eodory lemma, see
\cite{Sha1993}.
 We claim that $(\varphi_n)$ is a normal family
in $K$. Taking this for granted momentarily, we may assume that
$(\varphi_n)$ converges locally uniformly in $K$ to a holomorphic limit function
$\varphi : K \to \C$. It follows that
\begin{equation*}
 \varphi(z) =\frac{F(z)}{z F'(z)} \,  \quad \text{for } z \in \D \cap K \, ,  \end{equation*}
so $F(z)/(z F'(z))$ has an analytic continuation to $K$. In particular,
the zeros of the Blaschke product $F$ cannot accumulate on $I$, i.\!\;e.~$F$ has an analytic continuation across $I$, see
\cite[Chapter II, Theorem 6.1]{Gar2007}.

\medskip

It remains to prove that the sequence $(\varphi_n)$ is indeed a normal family
in $K$.
 Since $\varphi_n(\zeta) \in \R$ for any $\zeta \in I$, we have
\begin{equation*}
\varphi_n(\zeta)=\overline{\varphi_n(1/\overline{\zeta})} \quad \text{ for all
} \zeta \in K \, 
\end{equation*}
 by the Schwarz reflection principle. We now appeal to the  
Schwarz lemma  $$|F_n(z)| \le |z| \, , \qquad z \in \D \, , $$ and the Schwarz--Pick lemma
$$|F'_n(z)| \le \frac{1-|F_n(z)|^2}{1-|z|^2} \le \frac{1}{1-|z|^2} \, ,
\qquad z \in \D \, . $$
Using $\log^+x:=\max\{\log x, 0\}$ for $x \in \R$, $x >0$, these estimates lead to

\begin{equation*} 
 \begin{array}{rcl}
 \log^+ |\varphi_n(z)| & = & \log^+ \left| \displaystyle \frac{F_n(z)}{z F_n'(z)}  \right|
 \le 
\displaystyle \log ^+ \frac{1}{|F_n'(z)|}= \log^+ |F_n'(z)|-\log |F_n'(z)| \\[6mm]
& \le & \log \left(\displaystyle \frac{1}{1-|z|^2} \right)-\log |F_n'(z)|  \, , \quad z \in \D \cap K\, ,
\end{array}
\end{equation*}
for $ n=1,2,\ldots $. By reflection, we get the similar estimate
\begin{equation*} 
 \log^+ |\varphi_n(z)| \le \log^+ \frac{1}{|F_n'(1/\overline{z})|}
\, , \qquad z \in  K \backslash \D \, ,
\end{equation*}
for $ n=1,2,\ldots $.
We now choose a compact subset $\Omega$ of $K$. Then $\Omega \subset K_R(0)$
for some  $R>1$
and 
\begin{eqnarray*}
\iint \limits_{\Omega} \log^+  |\varphi_n(z)| \, d\sigma_z & \le &
\iint \limits_{\Omega \cap \D} \log^+  \frac{1}{|F_n'(z)|} \, d\sigma_z+\iint
\limits_{\Omega\backslash \D} \log^+
\frac{1}{|F'_n(1/\overline{z})|} \, d\sigma_z\\[1mm]
& \le & \iint \limits_{\D} \log^+  \frac{1}{|F_n'(z)|} \, d\sigma_z+\iint
\limits_{1 < |z|<R} \log^+
\frac{1}{|F'_n(1/\overline{z})|} \, d\sigma_z \\[1mm]
&\le&
 (1+R^4) \iint \limits_{\D} \log^+  \frac{1}{|F_n'(z)|} \, d\sigma_z \\[1mm]
& \le & (1+R^4) \iint \limits_{\D} \log \left(\frac{1}{1-|z|^2}\right) \, d\sigma_z
-(1+R^4) \iint \limits_{\D} \log |F_n'(z)| \, d\sigma_z \\[1mm]
& \le & c_1 -(1+R^4)  \iint \limits_{\D} \log |F_n'(z)| \, d\sigma_z 
\end{eqnarray*}

for some constant $c_1$ depending only on $\Omega$.
We note that $F_n'(z)=z^{N} g_n(z)$ with a holomorphic function
$g_n : \D \to \C$ that satisfies $g_n(0) \not=0$.
Consequently,
\begin{equation*} \iint \limits_{\Omega} \log^+  |\varphi_n(z)| \, d\sigma_z  \le  c_2-(1+R^4) \iint
\limits_{\D} \log |g_n(z)| \, d\sigma_z \, , \qquad n \ge 1 \, , \end{equation*}
where $c_2$ depends only on $\Omega$.
By the submean value inequality for subharmonic functions (see \cite[Theorem 2.6.8]{Ran}), we get
\begin{equation*} \iint \limits_{\Omega} \log^+  |\varphi_n(z)| \, d\sigma_z  \le  c_2-(1+R^4)  \pi 
 \log |g_n(0)|  \, , \qquad n \ge 1 \, . \end{equation*}
Since $g_n(0)=F^{(N+1)}_n(0)/(N+1)!$ and $F^{(N+1)}_n(0) \to
 F^{(N+1)}(0)>0$ as $n \to + \infty$ we obtain
\begin{equation*} \iint \limits_{\Omega} \log^+  |\varphi_n(z)| \, d\sigma_z  \le  c_3-(1+R^4) \pi 
 \log |F^{(N+1)}_n(0)| \le c_4< +\infty 
\end{equation*}
 for all $n\ge 1$, where $c_3$ and $c_4$ are some constants which do not depend
 on the functions $\varphi_n$.
Thus, if $\Omega' \subset K$ is compact and 
$$\text{dist}(z,\Omega'):=\inf\{ |z-z'| \, : \, z' \in \Omega'\}$$
denotes the euclidean distance of a point $z \in \C$ to $\Omega'$,
 then there is a $\delta>0$ such that the tubular neighborhood
 $\Omega_{\delta}:=\{z \in \C \, : \, \text{dist}(z,\Omega')\le \delta\}$ is
 entirely contained in $K$. Hence, there exists a constant $\tilde{c}$ such that for any
 $z \in \Omega'$ and $n=1,2, \ldots$,
\begin{equation*} \log^+|\varphi_n(z)| \le \frac{1}{\pi \delta^2} \iint \limits_{|w-z| \le \delta}
\log^+|\varphi_n(w)| \, d\sigma_w \le \frac{1}{\pi \delta^2} \iint
\limits_{\Omega_{\delta}}\log^+|\varphi_n(w)| \, d\sigma_w \le \frac{\tilde{c}}{\pi
  \delta^2} \, , \end{equation*}
where in the first inequality we used the submean value property of subharmonic
functions once more.
So $( \varphi_n)$ is  uniformly bounded on compact subsets of $K$ and
therefore  a normal family by
Montel's theorem. \hfill{$\blacksquare$}

\subsection{Proof of Theorem \ref{thm:union}}

Let $F$  and $G$ be maximal Blaschke
products with critical sets $A:=\mathcal{C}_1$ and $B:=\mathcal{C}_2$ respectively.
We fix a real number $c$ with $0<c<1$. Then
$$ \lambda_{A}(z) \, |dz|:=\frac{c \, |F'(z)|}{1-c^2 |F(z)|^2} \, |dz|\, , \qquad
\lambda_{B}(z) \, |dz|:=\frac{c \, |G'(z)|}{1-c^2 |G(z)|^2} \, |dz|$$
are two conformal pseudometrics on $\D$ with zero set $A$ resp.~$B$ and
constant curvature $-4$ on $\D \backslash A$ and $\D \backslash B$
respectively. The product pseudometric
$$ \lambda(z) \, |dz|:=\lambda_A(z) \lambda_B(z) \, |dz|$$
has zero set $A \cup B$ (counting multiplicities) and a computation shows that
\begin{equation} \label{eq:curvi1}
 \kappa_{\lambda}(z) =-4 \big[ \lambda_{A}(z)^{-2}+\lambda_B(z)^{-2} \big]
\, , \qquad z \in \D \backslash (A \cup B) \, .
\end{equation}  
We will show that there is a positive constant $\alpha$ such that
\begin{equation} \label{eq:hil10}
\kappa_{\lambda}(z) \le -\alpha \quad \text{ for all } z \in \D \backslash (A
\cup B)\, .
\end{equation}
Fix $ \xi \in \partial \D$. If $\xi \not \in \overline{A}$, then, by Theorem 
 \ref{thm:analytic_continuation}, the maximal Blaschke product
$F$ has an analytic continuation to a neighborhood $K_{\xi}$ of $\xi$.
Thus  $\lambda_{A}(z) \, |dz|$ is bounded on $K_{\xi} \cap \D$.
In a similar way, we see that the density $\lambda_B$ is bounded above near
any point $\xi \in \partial \D \backslash \overline{B}$. Since by hypothesis
$\overline{A} \cap \overline{B} \cap \partial \D=\emptyset$,  we conclude from
(\ref{eq:curvi1}) that the curvature $\kappa_{\lambda}$ is bounded above 
 close to any point on $\partial \D$. Since $\kappa_{\lambda}(z) \to -\infty$
 if $ z \to z_0 \in A \cup B$, there is positive constant $\alpha>0$ such that
(\ref{eq:hil10}) holds.

\medskip

 As a result of the estimate (\ref{eq:hil10}),
the conformal pseudometric
$$ \mu(z) \, |dz|:=\frac{\sqrt{\alpha}}{2} \, \lambda(z) \, |dz|$$
has curvature $\le -4$ on $\D \backslash (A \cup B)$ and zero set $A \cup B$.
Corollary \ref{cor:main0} shows that $A \cup B$ is an $H^{\infty}$ critical set.
\hfill{$\blacksquare$}
\subsection{Proof of Theorem \ref{thm:semigroup}}

Let $F,B$ be two maximal Blaschke products with critical sets $\mathcal{C}_F$
and $\mathcal{C}_B$, respectively. Consider 
the conformal pseudometric 
$$\sigma(z) \, |dz|:= ((B \circ F)^*\lambda_{\D})(z) \, |dz| \, ,$$
so
$$ \sigma(z) \, |dz|=\frac{|B'(F(z))| \, |F'(z)|}{1-|B(F(z))|^2} \, |dz| \, .$$
In order to show that $B \circ F$ is a maximal Blaschke product, we appeal to 
 Corollary \ref{cor:main1}, so 
we need to  show that $\sigma \equiv \sigma_{max}$, where $\sigma_{max}(z) \,
|dz|$ denotes the maximal conformal pseudometric on $\D$ with zero set
$\mathcal{C}_{B\circ F}$ and curvature $-4$ on $\D \backslash \mathcal{C}_{B
  \circ F}$. Recall that $\mathcal{C}_{B \circ F}$ denotes the critical set of
the function $B \circ F$.
 We proceed in two steps.

\medskip

{\it 1.~Step:} Assume $B$ is a {\it finite} Blaschke product.\\
Clearly, $\sigma \le \sigma_{max}$. 
In order to prove the reverse inequality we 
 define the auxiliary function
\begin{equation} \label{eq:w0}
 s(z):=\log \frac{\sigma_{max}(z)}{\sigma(z)} \, , \qquad z \in \D \, . 
\end{equation}
Since $\sigma_{max}(z) \, |dz|$ and $\sigma(z) \, |dz|$ do
have the same zero set $\mathcal{C}_{B \circ F}$, we first note that
$s : \D \to \R$ is continuous. Furthermore,  $$ \Delta \left( \log \frac{\sigma_{max}(z)}{\sigma(z)} \right)=
\Delta  \log \sigma_{max}(z)- \Delta  \log \sigma(z)
=4 \sigma_{max}(z)^2-4 \sigma(z)^2 \ge 0 $$
for every $z \in \D \backslash \mathcal{C}_{B \circ F}$. Since
$\mathcal{C}_{B\circ F}$ is a discrete subset of $\D$, we
see that the continuous function $s : \D \to \R$ is  subharmonic  on $\D$ and nonnegative by
construction. It remains to prove that $s \equiv 0$. 

\medskip

For this purpose we consider the following conformal pseudometrics,
\begin{align*}
\lambda_{max}(z) \, |dz| &:= F^*\lambda_{\D}(z) \,
|dz|=\frac{|F'(z)|}{1-|F(z)|^2} \,  |dz| \, \\
\intertext{and}
\mu_{max} (z) \, |dz| &:= B^*\lambda_{\D}(z) \,
|dz|=\frac{|B'(z)|}{1-|B(z)|^2} \, |dz| \, .
\end{align*}
By Corollary \ref{cor:main1},
 $\lambda_{max}(z) \, |dz|$  is the
maximal conformal pseudometric on $\D$ with zero set
$\mathcal{C}_F$ and curvature $-4$ on $\D \backslash \mathcal{C}_F$, and $\mu_{max} (z) \, |dz|$ is the
maximal conformal pseudometric on $\D$ with zero set
$\mathcal{C}_B$ and 
curvature $-4$ on $\D \backslash \mathcal{C}_B$.  

In particular, if $\tilde{B}$ denotes a finite
Blaschke product with {\it zero set} $\mathcal{C}_B$, then the conformal pseudometric
$$ |\tilde{B}(z)| \, \lambda_{\D}(z)\, |dz|=\frac{|\tilde{B}(z)|}{1-|z|^2} \, |dz|$$
has curvature $-4 |\tilde{B}(z)|^{-2} \le -4$ on $\D \backslash \mathcal{C}_B$, so we get the crucial estimate
$$ \mu_{max}(z) \ge  \frac{|\tilde{B}(z)|}{1-|z|^2} \qquad \text{ for any } z
\in \D \, .$$

\medskip

Since clearly, $\lambda_{max}(z) \ge \sigma_{max}(z) \ge \sigma(z)$ and since
$\sigma(z) \, |dz|$  can be written in the form
$$ \sigma(z) \, |dz|=\frac{\mu_{max}(F(z))}{\lambda_{\D}(F(z))} \cdot
\lambda_{max}(z) \, |dz| \, $$
the above estimate for $\mu_{max}(z)$ gives
$$ \lambda_{max}(z) \ge \sigma_{max}(z) \ge \sigma(z) \ge |\tilde{B}(F(z))|
\, \lambda_{max}(z) \, , $$
so rearranging terms yields
\begin{equation} \label{eq:w1}
 |\tilde{B}(F(z))| \le \frac{\sigma(z)}{\lambda_{max}(z)} \le
\frac{\sigma_{max}(z)}{\lambda_{max}(z)} \le 1 \, .
\end{equation}
We now make the obvious, but important observation that $\tilde{B} \circ F$
is a Blaschke product. This follows from the fact that $F$ as a maximal
Blaschke product is an {\it
  indestructible} Blaschke product and $\tilde{B}$ is a finite Blaschke
product. As a consequence of Frostman's characterization of Blaschke products (see \cite[Ch.~II, Theorem 2.4]{Gar2007}), we get
$$ \lim \limits_{r \to 1} \frac{1}{2 \pi} \int \limits_{0}^{2 \pi} \log |\tilde{B}(F(r e^{i
  t}))| \, dt=0 \, ,$$
which combined with (\ref{eq:w1}) leads to
$$  \frac{1}{2 \pi} \int \limits_{0}^{2 \pi} \log
\frac{\sigma_{max}(r e^{it})}{\sigma (r e^{it})} \, dt= \frac{1}{2 \pi} \int \limits_{0}^{2 \pi} \log
\frac{\sigma_{max}(r e^{it})}{\lambda_{max} (r e^{it})} \, dt - \frac{1}{2 \pi} \int \limits_{0}^{2 \pi} \log
\frac{\lambda_{max}(r e^{it})}{\sigma (r e^{it})} \, dt \to 0\, $$
as $r \to 1$. Recalling the definition of the nonnegative subharmonic function
$s : \D \to \C$ in (\ref{eq:w0}), we finally deduce that $s \equiv 0$.
This completes the proof that $B \circ F$ is a maximal Blaschke product for
the case that $B$ is a finite Blaschke product.

\bigskip

{\it 2.~Step:} Assume $B$ is not a finite Blaschke product.\\
Let $\mathcal{C}_B=(z_j)$. Then for each positive integer $n$ there exists a
finite Blaschke product $B_n$ with critical set $\mathcal{C}_{B_n}=(z_1, \ldots,
z_n)$, see Remark \ref{rem:wu} (a).
From the first step, we know that $B_n \circ F$ is a maximal Blaschke
product with critical set $\mathcal{C}_{B_n \circ F} \subset \mathcal{C}_{B \circ F}$, so
$$ \frac{|(B_n \circ F)'(z)|}{1-|(B_n \circ F)(z)|^2} \ge \sigma_{max}(z) \ge 
\frac{|(B \circ
  F)'(z)|}{1-|(B \circ F)(z)|^2} \, , \qquad z \in \D \, .$$
Since Corollary \ref{cor:main} shows that the finite Blaschke products $B_n$
converge locally uniformly in $\D$ to $B$,  we get
$$ \sigma_{max}(z) =
\frac{|(B \circ
  F)'(z)|}{1-|(B \circ F)(z)|^2} \, , \qquad z \in \D \, , $$
as desired. \hfill{$\blacksquare$}

\subsection{Proof of Theorem \ref{thm:bdd}}
Let $B$ be a maximal Blaschke product with critical set $\mathcal{C}$
and let $\zeta \in \partial \D \backslash \overline{\mathcal{C}}$. 
By Theorem \ref{thm:analytic_continuation}, the function $B$ has an
analytic continuation across an open arc $\Gamma \subseteq \partial \D$ which
contains $\zeta$ and such that 
$B(\Gamma) \subset \partial \D$. Hence, the so--called boundary Ahlfors' lemma (see
Theorem 1.1 in \cite{KRR06}) shows that
condition (\ref{eq:bdd}) holds for the boundary point $\zeta$.

\section{Extension to more general regions}

In this section we extend Theorem \ref{thm:2} and 
Theorem \ref{thm:analytic_continuation} to more general subdomains of the
complex plane than the unit disk. 
For this purpose we 
note that  the extremal problem (\ref{eq:ex1}) of Theorem \ref{thm:2} 
has a natural counterpart for 
functions analytic and bounded in other domains than the unit disk:

\begin{remark}[The Riemann mapping theorem and the Ahlfors mapping theorem] \label{rem:2}
Let $\Omega \subseteq \C$ be a domain containing the origin and let $\mathcal{C}=(z_j)$
  be the  critical set of a nonconstant function $f$ in $H^{\infty}(\Omega)$. Here,
   $H^{\infty}(\Omega)$ is  the set of all functions analytic and
  bounded in $\Omega$. We denote by $N$  the number 
of times that $0$ appears in the sequence ${\cal  C}$ and let 
$$ \mathcal{F}_{\mathcal{C}}(\Omega):=\left\{f \in H^{\infty}(\Omega) \, : \, f'(z)=0 \text{ for
  any } z \in \mathcal{C}\right\} \, .$$
Then, by a normal family argument, there is always at least one extremal
function for the extremal problem
\begin{equation} \label{eq:ex2}
 \sup \big\{ \Re f^{(N+1)}(0) \, : f \in {\cal F}_{\cal C}(\Omega), \, {||f||}_{\infty} \le 1 \big \} \, .
\end{equation}
For the special case $\mathcal{C}=\emptyset$, 
 this extremal problem  leads to the
Riemann resp.~Ahlfors mapping function:
\begin{itemize}
\item[(a)] {\it $\mathcal{C}=\emptyset$ and $\Omega \subsetneq \C$  simply
    connected }\\
 In this case  the 
  extremal problem (\ref{eq:ex2}) has a unique extremal function $\Psi$,
  namely the normalized Riemann map
  $\Psi$ of the domain $\Omega$, which maps $\Omega$ conformally onto $\D$ such  that 
  $\Psi(0)=0$ and $\Psi'(0)>0$. In order to see this, we note that the standard textbook proof of existence of
  the Riemann map is through
  maximzing $\Re f'(0)$ over all  {\it injective} functions in $\mathcal{F}_{\emptyset}(\Omega)$.
  That this apparently more restrictive extremal problem has the same solution
  as the original problem follows from the Schwarz lemma. To verify this, we first
  note that it is easily seen that the extremal function $\phi : \Omega \to \D$ to (\ref{eq:ex2}) is
  normalized by $\phi(0)=0$ and $\phi'(0)>0$. Since the normalized Riemann map
  $\Psi$ of $\Omega$ belongs to $\mathcal{F}_{\emptyset}(\Omega)$, we get
  $\Psi'(0) \le \phi'(0)$, so the holomorphic self--map $\omega:=\Psi^{-1} \circ \phi$
  of the unit disk fixes the origin and satisfies $\omega'(0) \ge 1$. 
 The Schwarz lemma now implies that $\Psi^{-1} \circ \phi : \D \to \D$ is the
 identity function, so $\phi=\Psi$. We don't know of a {\it direct} way of
 proving that the extremal function $\phi$ for the extremal problem
 (\ref{eq:ex2}) is a conformal map from the domain $\Omega$ onto the unit disk
 $\D$.

\item[(b)] {\it $\mathcal{C}=\emptyset$ and $\Omega \subsetneq \C$ a smooth
    multiply connected domain}\\ 
If the domain $\Omega$ has finite connectivity $n \ge 2$, none of whose boundary components
reduces to a point, then the
  extremal problem (\ref{eq:ex2}) has a unique extremal function, namely the Ahlfors map
  $\Psi : \Omega \to \D$. It is an $n:1$ map  from $\Omega$ onto $\D$ such  that 
  $\Psi(0)=0<\Psi'(0)$; see \cite{Ahl}, \cite{Gru1978} and \cite{Bell1992}. We also refer to \cite{BeuAhl} where
  the same extremal problem was considered in generality for various classes
  of analytic functions.
\end{itemize}
\end{remark}

When $\Omega$ is a simply connected domain that is not equal to the whole
complex plane, then Theorem \ref{thm:2}
combined with the Riemann mapping theorem (Remark \ref{rem:2} (a))  leads to 
the following result about the extremal problem (\ref{eq:ex2}).

\begin{theorem} \label{cor:thm1}
Let $\Omega \subsetneq \C$ be a simply connected domain containing the origin, let
$\mathcal{C}=(z_j)$ be an
$H^{\infty}(\Omega)$ critical set and let $N$ denote the number of times that
$0$ appears in the sequence $\mathcal{C}$. Then the extremal problem
(\ref{eq:ex2}) has the unique extremal function $B_{\Psi(\mathcal{C})} \circ
\Psi$. Here, $\Psi$ is the normalized Riemann map for the domain $\Omega$ and
$B_{\Psi(\mathcal{C})}$ is the extremal function for the extremal problem (\ref{eq:ex1})
for the critical set $\Psi(\mathcal{C})$ 
according to Theorem \ref{thm:2}.
\end{theorem}

\hide{If $\mathcal{C}\not=\emptyset$ and
 $\Omega \subsetneq \C$ is simply connected, then the remarks in (a) and Theorem
 \ref{thm:2} give the following information.
If $\Psi$ is the normalized Riemann map for the domain $\Omega$ and
$B_{\Psi(\mathcal{C})}$ is the extremal function for the extremal problem (\ref{eq:ex1})
for the critical set $\Psi(\mathcal{C})$ 
according to Theorem \ref{thm:2}, then, clearly,
$B_{\Psi(\mathcal{C})} \circ \Psi$ is the unique extremal function for the
extremal problem (\ref{eq:ex2}).
Hence the extremal function the critical set $\mathcal{C}$ is the composition of a maximal Blaschke for 
$\Psi(\mathcal{C})$ and the conformal map $\Psi$.}

It would be interesting to study the extremal problem (\ref{eq:ex2}) in the
presence of critical points  when
  the domain $\Omega$ has connectivity $n>1$, none of whose boundary components
reduces to a point.

\medskip

The next result is an immediate consequence of Theorem \ref{thm:2} and the Schwarz
reflection principle.

\begin{theorem}
Let $\Omega \subsetneq \C$ be a simply connected domain containing the origin
and bounded by a 
real--analytic Jordan arc, and let $\mathcal{C}$ be an $H^{\infty}(\Omega)$ critical
set. Then the extremal function for the extremal problem (\ref{eq:ex2}) has an
analytic continuation across any point of $\partial \Omega \backslash \overline{\mathcal{C}}$.
\end{theorem}

Perhaps the same result holds for any extremal function for the extremal
problem (\ref{eq:ex2}) when
 $\Omega$ is a multiply connected domain with connectivity $n>1$ and bounded by $n$
real--analytic non--intersecting Jordan curves.\medskip

\section{Further remarks and open problems}

We close this paper with a number of remarks and further open problems.

\medskip

Let us first return to Theorem \ref{thm:2}.
It says that every maximal Blaschke product is
indestructible. This unvoidably  suggests the following question.

\begin{problem} \label{prob:1}
Is every indestructible Blaschke product a maximal Blaschke product\,?
\end{problem}

Note that an affirmative answer to Problem \ref{prob:1} would in particular imply that
 every locally univalent indestructible Blaschke product is in fact  univalent.

\medskip

The semigroup property of  maximal Blaschke products (Theorem
\ref{thm:semigroup}) states that if $B,C \in H^{\infty}$ are maximal Blaschke
products then the composition
$A:=B \circ C$ is  again a maximal Blaschke product. 
It is natural to ask for the converse statement: If $B,C \in H^{\infty}$ and
$A=B \circ C$ is a maximal
Blaschke product, does it follow that $B$ and $C$ are
also maximal Blaschke products\,?
The following simple observation provides a partial answer.

\begin{proposition} \label{prop:semi1}
Let $A,B,C \in H^{\infty}$ such that $A=B \circ C$ is 
a maximal Blaschke
product.
Then $B$ is a maximal Blaschke product.
\end{proposition}

\textit{Proof.} 
Let $\tilde{B}$ be a maximal Blaschke product with critical set 
$\mathcal{C}_B$ and let $\lambda_{max}(z) \, |dz|$
denote the maximal conformal pseudometric on $\D$ with zero set
$\mathcal{C}_{B}$ and curvature $-4$ on $\D\backslash \mathcal{C}_B$. Then, by
Corollary \ref{cor:main1},
$$ \lambda_{max}(w)=\frac{|\tilde{B}'(w)|}{1-|\tilde{B}(w)|^2} \ge \frac{|B'(w)|}{1-|B(w)|^2}
\quad \text{ for every } w \in \D \, .$$
This inequality  holds in particular for any $w=C(z)$, $z \in \D$, 
so  by multiplying both sides with $|C'(z)|$, we get for every $z \in \D$
$$ C^*\lambda_{max}(z)=\frac{|(\tilde{B} \circ C)'(z)|}{1-|(\tilde{B} \circ C)(z)|^2} \ge
\frac{|(B \circ C)'(z)|}{1-|(B \circ C)(z)|^2}=\frac{|A'(z)|}{1-|A(z)|^2} \ge C^*\lambda_{max}(z)\, ,$$
where the last inequality comes from the maximality of $A$. Hence  equality
holds throughout and 
Liouville's
Theorem implies  $B \circ C=T \circ \tilde{B} \circ C$ for some unit disk automorphism
$T$, i.e., $B=T \circ \tilde{B}$. Thus $B$ is a maximal Blaschke product.
\hfill{$\square$}

\begin{problem} \label{prob:2}
Let $A,B,C \in H^{\infty}$ such that $A=B \circ C$ is 
a maximal Blaschke product. Does it follow that $C$ is a maximal Blaschke product\,?
\end{problem}

\begin{remark}
It is well--known and easy to prove that if the function
$A=B \circ C$
in Problem \ref{prob:2} is a \textit{finite} Blaschke product, then  
both $B$ and $C$ are finite Blaschke products.
Hence, if the answer to Problem \ref{prob:2} is affirmative, then 
it would be interesting to explore the possibility of a ``prime factorization''
of maximal Blaschke products in 
a way similar to the recent extension of
Ritt's celebrated factorization theorem for finite
Blaschke products due to Ng and Wang (see \cite{NgWang2011}). 
In this connection, information about the critical \textit{values} of
maximal Blaschke products would be valuable.
\end{remark}

\begin{remark} 
In view of Problem \ref{prob:1} and the discussion above, the following two
questions arise.
\begin{itemize}
\item[(a)]
Is the composition of two indestructible Blaschke products $B$ and $C$ an
indestructible Blaschke product\,?
\item[(b)] If $B,C \in H^{\infty}$ and $B \circ C$ is an indestructible
Blaschke product, must also $B$ and $C$ be indestructible Blaschke products\,?
\end{itemize}

These issues will be discussed in the forth--coming paper \cite{KR}.
\end{remark}

\vfill
\hspace{0.9cm}\begin{minipage}{8cm}
Daniela Kraus and Oliver Roth\\
Department of Mathematics\\
University of W\"urzburg\\
97074 W\"urzburg\\
Germany\\
dakraus@mathematik.uni-wuerzburg.de\\
roth@mathematik.uni-wuerzburg.de

\end{minipage}

\end{document}